\documentclass[11pt]{amsart}
\usepackage{txfonts}
\usepackage{}
\usepackage{amssymb}
\theoremstyle{plain}
\newtheorem{Thm}{Theorem}

\begin{document}

\title[bounded states of one dimensional Schrodinger systems]
{bounded states of one dimensional Schrodinger systems}

\author{Li Ma}
\address{Distinguished Professor, Department of mathematics \\
Henan Normal university \\
Xinxiang, 453007 \\
China} \email{lma@tsinghua.edu.cn}

\thanks{The research is partially supported by the National Natural Science
Foundation of China 10631020 and SRFDP 20090002110019}

\begin{abstract}
In this paper, we study the existence problem of bound states of one
dimensional Schrodinger system via the blow-up method.

{ \textbf{Mathematics Subject Classification 2000}: 53Cxx,35Jxx}

{ \textbf{Keywords}: Schrodinger system, stationary solution, one
dimensional}
\end{abstract}

 \maketitle

\section{Introduction}
The coupled nonlinear Schrodinger (NLS) equations arise naturally in
nonlinear optics and they are also considered as the model of
Bose-Einstein condensate. As a interesting physical model, the
planar stationary light beams propagating in the z-direction in a
nonlinear medium are described by an NLS equation of the form
\begin{equation}\label{NSL}
iE_z+E_{xx}+k|E|^2E=0
\end{equation}
where $i$ denotes the imaginary unit, $E=E(x,z)$ denotes the complex
envelope of an electric field and $E_{xx}$ is its second partial
derivative, and $k$ is the coupling constant, which is assumed to be
positive and corresponds physically self-focusing of the medium. If
we are looking for the bounded states (or standing waves), namely
the special solutions to the system above in the form
$$
E_j(z,x)=e^{i\lambda z}u_j(x), \ \ \lambda>0 ,
$$
we are led to solve the elliptic system of the form
$$
-u_{xx}+\lambda u=k|u|^2u, \ \ x\in \mathbf{R}
$$
where $u=(u^1,u^2)$ is the unknown vector valued function on the
real line $\mathbf{R}$. The cubic nonlinear term on the right hand
side makes the system possible to be solved explicitly. When the
physical situation is more complicated, the coupling constant $k$
can be changed into a positive function $Q(t)$ and we are led to
consider the elliptic system
$$
-u_{xx}+\lambda u=Q(t)|u|^2u, \ \ \ x\in \mathbf{R}.
$$

Let $Q(t)=Q(|t|)\geq 0$ be a bounded even function on $\mathbf{R}$.
Let $a>0$ be a positive constant. Let $p>1$. Assume that $Q(0)>0$.
We shall consider the non-trivial non-negative solutions to the
following elliptic system
\begin{equation}\label{sch}
-u_{tt}+a u=Q(t)|u|^{p-1}u, \ \ \ x\in \mathbf{R}.
\end{equation}
The study of the problem above is non-trivial in the sense that
there is no Sobolev imbedding theorem on the whole real line. This
fact is overlooked by some previous authors.

There are some results concerning with the special case when $p=3$
by the perturbation method, see \cite{A}. The higher dimensional
cases have studied by Lin-Wei \cite{LW}, Ma-Zhao \cite{MZ} and
others \cite{BWW}, \cite{Ma},\cite{LM}, \cite{W}.

We set up the following result.
\begin{Thm}\label{mali} There is at least one non-negative non-trivial bounded solution
$u$ to (\ref{sch}).
\end{Thm}

The proof of this result is by the blow-up method. Roughly speaking,
it is a proof of arguing by contradiction, which is given in next
section. A general result can be proved similarly and is stated in
the last section.

\section{the proof of Theorem \ref{mali}}

For any $R>0$, we denote $I_R=[-R,R]$. Define on $H_R=H_0^1(I_R)$,
the functional
$$
E_R(u)=\int_{I_R} \frac{1}{2}(|u'|^2+au^2)-\frac{1}{p+1} \int_{I_R}
Q(t)|u|^{p+1}.
$$
Its Euler-Lagrange equation is
$$
-u^{''}+au=Q(t)|u|^{p-1}u, \ \ in \ \  I_R
$$
with $u=0$ at $t=R$ and $-R$.

Introduce the Nehari manifold
$$
\mathbf{N}=\{u\in H_R-\{0\}; \int_{I_R}(|u'|^2+au^2)=\int_{I_R}
Q(t)|u|^{p+1}\}.
$$

Restricted to $\mathbf{N}$ the functional $E_R$ can be written as
$$
E_R(u)=(\frac{1}{2}-\frac{1}{p+1})\int_{I_R}(|u'|^2+au^2).
$$

Then by the direct method we know that there is a minimizer $u_R\in
\mathbf{N}$ to the minimization problem
$$
d_R:=\int_{\mathbf{N}}E_R(u)>0.
$$
Furthermore, $u_R>0$ in the sense that each component is positive.
In fact using the Sobolev inequality we know that
$$
\int_{I_R} Q(t)|u|^{p+1}\leq C[\int_{I_R}(|u'|^2+au^2)]^{(p+1)/2}.
$$

For $u\in \mathbf{N}$, we have
$$
\int_{I_R}(|u'|^2+au^2)=\int_{I_R} Q(t)|u|^{p+1}.
$$
Then we have
$$
\int_{I_R}(|u'|^2+au^2)\leq C[\int_{I_R}(|u'|^2+au^2)]^{(p+1)/2}.
$$
This implies that
$$
\int_{I_R}(|u'|^2+au^2)\geq C(R)>0
$$
for some uniform constant $C(R)>0$. Hence $d_R\geq C(R)>0$. We may
assume that each component is non-trivial, otherwise, we just delete
the trivial functions.

Using the symmetry property we know that each component of $u_R$ has
its maximum point at origin. For otherwise the component will be
trivial and then it is zero.

Without loss of generality we assume that the maximum component of
$u_R$ is the first one in the sense that
$$
u^1_R(0)=\max_{I_R} u_R^1(t)=M_R\geq u^j_R(0)>0.
$$

 Let $R=R_j\to\infty$. Let $M_j=M_{R_j}$ and $u_j=u_{R_j}$. We \emph{claim} that $M_j$ is uniformly bounded
 and then we have proved Theorem \ref{mali} by taking the limit.
 Assume that $M_j\to \infty$, which will be showed that it is impossible.
We define the blow-up sequence
 $$
v_j(t)=M_j^{-1}u_j(M_j^{\beta}t)
 $$
 for $\beta>0$ to be determined later. Let
 $$
Q_j(t)=Q(M_j^{\beta}t).
 $$

Note that
$$
v_j^{''}(t)=M_j^{2\beta-1}u_j^{''}
$$
and
$$
|v_j|^{p-1}v_j=M_j^{-p}|u_j|^{p-1}u_j.
$$
We let $\beta=\frac{1-p}{2}<0$. Since $v_j$ is uniformly bounded, we
may use the elliptic regularity theory to conclude that $(v_j)$ has
a convergent subsequence $(v_{j_k})$ in $C^3_{loc}$, which is still
denoted by $(v_j)$ and its  $C^2$ limit $\bar{v}$ such that
\begin{equation}\label{limit}
-\bar{v}^{''}=Q(0)\bar{v}^p, \ \ in \ \ (-\infty, \infty)
\end{equation}
with $\bar{v}^1(0)=1$.

We shall prove that $(\ref{limit})$ has no nontrivial positive
solution. Let $f=\bar{v}^1$. by maximum principle we know that $f>0$
in $(-\infty,\infty)$. \emph{Claim} that $f'>0$. If not, then for
some $t_0$, we have
$$
f'(t_0)\leq 0.
$$
Then we have for $t>t_0$,
$$
f'(t)=f'(t_0)+f^{''}(\bar{t})(t-t_0)<0.
$$
Here we have used $f^{''}(\bar{t})=Q(0)|\bar{v}|^{p-1}f<0$.

For any $t>t_1>t_0$, we have
$$
f'(t)<f'(t_1)<0.
$$
Hence we have
$$
f(t)=f(t_1)+\int_{t_1}^t f'(s)ds\leq
f(t_1)+f'(t_1)(t-t_1)\to-\infty.
$$
as $t\to \infty$, a contradiction.

We remark that the same argument works for other component. Then we
have $|\bar{v}|^{p-1}f(t)\geq |\bar{v}|^{p-1}f(s)$ for $t>s$.

Then,
$$
-f^{''}(t)=Q(0)|\bar{v}|^{p-1}f(t)\geq Q(0)|\bar{v}|^{p-1}f(s)>0.
$$

We then conclude that
$$
f'(s)>f'(s)-f'(t)=Q(0)\int_s^t|\bar{v}|^{p-1}f(\tau)d\tau\geq
Q(0)|\bar{v}|^{p-1}f(s)(t-s)\to\infty
$$
as $t\to\infty$, which is impossible.

\section{Remarks about related systems}
One may also look for bounded states of the form
$$
E_j(z,x)=e^{i\lambda_j z}u_j(x), \ \ \lambda_j>0,
$$
 to the system
(\ref{NSL}) and the reduced system is of the form
\begin{equation}\label{modify}
-u^j_{xx}+\lambda_j u^j=Q(t)|u|^{p-1}u^j, \ \ \ x\in \mathbf{R}. \ \
\end{equation}
Here $p>1$ and $\lambda_j>0, \ \ j=1,...,N$.

Then using similar method to theorem \ref{mali}, we may use the
functional
$$E_R(u)=\int_{I_R}
\frac{1}{2}(|u'|^2+\sum_j\lambda_j(u^j)^2)-\frac{1}{p+1} \int_{I_R}
Q(t)|u|^{p+1}
$$
to obtain the following result.

\begin{Thm}\label{mali2} Assume that $Q(t)=Q(|t|)\geq 0$ is a bounded function in $\mathbf{R}$ with $Q(0)>0$
and $\lambda_j>0$ for all $j=1,...,N$. There is at least one
non-negative non-trivial bounded solution $u=(u^1,...,u^N)$ to the
system (\ref{modify}).
\end{Thm}

\end{document}